\documentclass[a4paper,10pt,twoside]{amsart}
\usepackage[utf8]{inputenc}
\usepackage{latexsym}
\usepackage[pdftex]{graphicx}
\usepackage{hyperref}
\usepackage{amsmath,amssymb, amsthm, amscd}
\usepackage{dsfont}
\usepackage{multicol}
\usepackage[marginratio=1:1,height=597pt,width=460pt,tmargin=117pt]{geometry}
\usepackage[all]{xy}
\usepackage{etoolbox}
\usepackage{enumitem}
\usepackage{tikz-cd}
\usepackage{mathtools}
\usepackage{float}
\usepackage{comment}
\usepackage{array}
\usepackage{bm}
\usepackage[pagewise]{lineno}

\expandafter\patchcmd\csname\string\proof\endcsname
  {\normalparindent}{0pt }{}{}
\makeatletter
\patchcmd{\@thm}{\thm@headfont{\scshape}}{\thm@headfont{\scshape\bfseries}}{}{}
\patchcmd{\@thm}{\thm@notefont{\fontseries\mddefault\upshape}}{}{}{}
\makeatother
\makeatletter
\patchcmd\@thm
  {\let\thm@indent\indent}{\let\thm@indent\indent}%
  {}{}
\makeatother
\newtheorem*{theorem*}{Theorem}
\newtheorem{theorem}[equation]{Theorem}
\newtheorem{lemma}[equation]{Lemma}
\newtheorem{proposition}[equation]{Proposition}

\theoremstyle{definition}

\theoremstyle{definition}
\newtheorem{remark}[equation]{Remark}
\theoremstyle{remark}

\numberwithin{subsection}{section}
\numberwithin{equation}{section}
\newcommand{\iso}{\xrightarrow{
   \,\smash{\raisebox{-0.50ex}{\ensuremath{\scriptstyle\sim}}}\,}}

\newcommand{\tld}{\widetilde}
\newcommand{\cO}{\mathcal{O}}
\newcommand{\Z}{\mathbb{Z}}
\newcommand{\F}{\mathbb{F}}
\newcommand{\Q}{\mathbb{Q}}
\newcommand{\ovl}{\overline}
\DeclareMathOperator{\GL}{\mathrm{GL}}
\DeclareMathOperator{\Hom}{\mathrm{Hom}}
\DeclareMathOperator{\Ind}{\mathrm{Ind}}
\DeclareMathOperator{\End}{\mathrm{End}}
\DeclareMathOperator{\soc}{\mathrm{soc}}
\DeclareMathOperator{\cosoc}{\mathrm{cosoc}}
   
\title{Irreducible smooth representations in defining characteristic without central character}
\author{Daniel Le}
\address{Department of Mathematics,
Purdue University,
150 N. University Street, 
West Lafayette, IN 47907-2067}
\email{ledt@purdue.edu}

\begin{document}
\maketitle
\begin{abstract}
Let $p>3$, $n>1$ an integer, and $F$ be a non-archimedean local field with residue field a proper finite extension of $\mathbb{F}_p$. 
Let $E$ be an algebraically closed countable field extension of the residue field of $F$. 
In this short note, we explain how the methods from \cite{le,GLS} can be used to construct irreducible smooth representations of $\mathrm{GL}_n(F)$ over $E$ without a central character. 
We also construct irreducible smooth representations of $\mathrm{GL}_n(F)$ over $E$ with simultaneously a central character, nonscalar endomorphisms, and if $n>3$, without a Hecke eigenvalue. 

\end{abstract}

\section{Introduction}

Let $p$ be a prime number. 
The study of smooth representations of $p$-adic reductive groups grew out of the theory of automorphic forms. 
An important property of smooth representations that naturally arises in that theory is \emph{admissibility}: the space of invariants for any compact open subgroup is finite-dimensional. 
Admissibility implies a number of other desirable properties. 
Let $F$ be a nonarchimedean local field with ring of integers $\cO_F$ and residue characteristic $p$, and let $G$ be the group of $F$-points of a reductive group. 
Fix an algebraically closed field $E$ and an irreducible smooth representation $\pi$ of $G$ over $E$. 
For brevity, let us say that $\pi$ is \emph{central} if it has a central character, and that it is \emph{Schur} if the $E[G]$-endomorphisms of $\pi$ are scalar i.e.~in $E$. 
It is not hard to see that $\pi$ is admissible implies that $\pi$ is Schur and that $\pi$ is Schur implies that $\pi$ is central. 
If the characteristic of $E$ is distinct from $p$, then $\pi$ is irreducible implies that $\pi$ is admissible \cite[II.2.8]{Vigneras} and thus Schur and central. 
The proof relies crucially on the existence of $E$-valued Haar measures which do not exist when the characteristic of $E$ is $p$. 
The characteristic $p$ case has received increasing attention because of its relation to mod $p$ automorphic forms and the mod $p$ Langlands program. 
We assume from now on that the characteristic of $E$ is $p$. 

When $G = \GL_2(\Q_p)$, Barthel, Livn\'e, and Breuil \cite{BL,Breuil} showed that admissible, Schur, and central are all equivalent. 
A crucial tool in their classification is \emph{Hecke theory}. 
Let $n$ be a positive integer, and let $G = \GL_n(F)$. 
(Much of the discussion below can be extended to general reductive groups, but we will restrict to $\GL_n$ for simplicity and also because that is the context for the rest of the article.) 
A (\emph{Serre}) \emph{weight} is an irreducible smooth $E[\GL_n(\cO_F)]$-representation. 
We say that $\pi$ admits a Hecke eigenvalue (or is \emph{Hecke}) if $\pi$ is isomorphic to a quotient of a representation of the form 
$
(\mathrm{ind}_{\GL_n(\cO_F)}^{\GL_n(F)} \tau) \otimes_{\mathcal{H}(\tau)} \chi
$
where $\tau$ is a weight, $\mathrm{ind}$ denotes compact induction, $\mathcal{H}(\tau)$ is the Hecke algebra $\End_{E[G]}(\mathrm{ind}_{\GL_n(\cO_F)}^{\GL_n(F)}\tau)$, and $\chi$ is an $E$-homomorphism $\mathcal{H}(\tau) \rightarrow E$. 
One can show that $\pi$ is admissible implies that $\pi$ is Hecke (using that $\mathcal{H}(\tau)$ is commutative and finitely generated over $E$; see \cite{herzig-satake}) and $\pi$ is Hecke implies that $\pi$ is central (since $(\mathrm{ind}_{\GL_n(\cO_F)}^{\GL_n(F)} \tau) \otimes_{\mathcal{H}(\tau)} \chi$ has a central character). 
When $n=2$, Barthel--Livn\'e showed further that $\pi$ is central implies that $\pi$ is Hecke. 
We expect that similar methods would yield the analogous result for groups of semisimple rank at most $1$ (see e.g.~\cite{xu}). 
For $\GL_2(\Q_p)$, Breuil showed that $\pi$ is Hecke implies that $\pi$ is admissible. 
Thus, for irreducible smooth representations of $\GL_2(\Q_p)$, admissible, Schur, Hecke, and central are all equivalent. 
Later, Berger \cite{Berger}, building on these results, showed that all irreducible smooth representations of $\GL_2(\Q_p)$ are central so that all of these properties hold. See Figure \ref{figure:imply} for a summary of the above results.

\begin{figure}\label{figure:imply}
\begin{center}
\begin{tikzcd}[arrows=Rightarrow]
& \textrm{Schur} \arrow[dr]
& & \\
\textrm{admissible} \ar[ur] \ar[dr]
&
& \textrm{central} \ar[dl, bend left,"\GL_2(F)"] & \textrm{irreducible} \ar[l,swap,"\GL_2(\Q_p)"]\\
& \textrm{Hecke} \ar[ur] \ar[ul,bend left,"\GL_2(\Q_p)"]
&
\end{tikzcd}
\caption{Some known implications for irreducible smooth representations of $p$-adic reductive groups over a characteristic $p$ algebraically closed field. }
\end{center}
\end{figure}

It is natural to ask whether all irreducible smooth $E$-representations of the group of $F$-points of a reductive group are admissible, Schur, Hecke, or central \cite[Questions 1, 2, and 8]{ahhv}, or even the extent to which any of these properties implies any others. 
This is of particular interest because Abe, Henniart, Herzig, and Vign\'eras \cite{herzig,AHHV-JAMS} have generalized the Barthel--Livn\'e classification for $\GL_2$ to any reductive group. 
In this article, we show that the implications in Figure \ref{figure:imply} are essentially the best possible. 
To this end, we construct examples, building on methods of Ghate, Le, and Sheth \cite{le,gs,GLS}, of mod $p$ irreducible smooth representations of $\GL_n(F)$ with a mix of different properties. 

\begin{theorem} \label{thm:main}
    Let $n>1$ and $p>3$ be a prime. 
    Let $F$ be a nonarchimedean local field with residue field a proper extension of $\F_p$. 
    \begin{enumerate}
        \item \label{item:central}  
        Let $E$ be an extension of the residue field of $F$ of countable cardinality.
        Then there exists an irreducible smooth $\GL_n(F)$-representation over $E$ which does not have a central character. 
        \item \label{item:endo}  
        Let $E$ be an extension of the residue field of $F$ and $\tld{E}$ be a division algebra over $E$ of countable dimension with center $\tld{C}$. 
        Then there exists an irreducible nonadmissible smooth $\GL_n(F)$-representation over $E$ whose endomorphism algebra is isomorphic to $\tld{E}$ as $E$-algebras. 
        Moreover, there exists such a representation
        \begin{itemize}
            \item with a Hecke eigenvalue (and a central character),
            \item without a central character (and without a Hecke eigenvalue) if $E\subsetneq \tld{C}$, or
            \item with a central character but without a Hecke eigenvalue if $n>3$ and $E\subsetneq \tld{C}$. 
        \end{itemize}
    \end{enumerate}
\end{theorem}

\begin{remark}\mbox{}
    \begin{enumerate}
        \item Take $E$ to be algebraically closed and countable. 
        Under the hypotheses of Theorem \ref{thm:main}, \eqref{item:central} shows that for smooth $\GL_n(F)$-representations over $E$, irreducible does not imply central. 
        Further, \eqref{item:endo} shows that central does not imply Schur (take $\tld{E}$ to be a proper countable field extension of $E$, e.g.~the field $E(X)$ of rational functions in one variable), central does not imply Hecke (when $n>3$; e.g.~take $\tld{E} = \tld{C}$ to be $E(X)$ again), Hecke does not imply admissible or Schur (again take $E \subsetneq \tld{E}$), and Schur does not imply admissible (take $\tld{E} = E$). 
        We do not know whether Schur (or just that the center of the endomorphism algebra is $E$) implies Hecke except that it does when $n=2$. 
        Given the results for $\GL_2(\Q_p)$, the hypothesis that the residue field of $F$ is a proper extension of $\F_p$ cannot be completely removed, though we expect that examples with similar properties can also be constructed for $\GL_n(F)$ with $F\not\cong \Q_p$. 
        \item If $\tld{E}$ has infinite dimension over $E$, then nonadmissibility in Theorem \ref{thm:main}\eqref{item:endo} follows automatically from the endomorphism statement and the injectivity of $\End_{E[\GL_n(F)]}(\pi) \rightarrow \End_{E}(\pi^K)$ for a compact open pro-$p$ subgroup $K$. 
        \item Taking $E$ to be the residue field of $F$ and $\tld{E}$ the algebraic closure of $E$ in Theorem \ref{thm:main}\eqref{item:endo} recovers \cite[Theorem 1.2]{GLS}. 
        \item The endomorphism algebra of $\pi$ is necessarily a division algebra over $E$ of countable dimension, since $\pi$ is (generated by any nonzero vector and) itself of countable dimension over $E$ being the quotient of a compact induction from a compact open subgroup of a finite-dimensional representation. 
        \item If $E$ is algebraically closed and \emph{uncountable}, then every irreducible smooth representation of a $p$-adic reductive group satisfies Schur's lemma and thus has a central character. 
        Indeed, a non-scalar endomorphism $\varphi$ of such a representation $\pi$ would give an injection of the function field $E(\varphi)$ in one variable into the endomorphism algebra of $\pi$ which, as before, is of countable dimension over $E$. 
        However, $E(\varphi)$ has uncountable dimension over $E$ since $\{\frac{1}{\varphi-\alpha} \mid \alpha \in E\}$ is a linearly independent set. 
        We then know in this case that $\pi$ is Hecke if $n=2$, but we do not know this if $n>2$.  
    \end{enumerate}
\end{remark}

\subsection{Acknowledgments}
Some key ideas in this note appeared in \cite{le}, and we embarrassingly did not notice the applications to questions about central characters at that time. 
We thank Benjamin Schraen again for asking about endomorphism algebras years ago. 
We thank Eknath Ghate and Mihir Sheth for comments on previous drafts of this paper and for our earlier collaboration. 
We thank Florian Herzig for the suggestion to use ordinary parts in the proof of Theorem \ref{thm:endo}. 
Finally, we thank the anonymous referee for many comments and corrections, for the suggestion to include a computation of the endomorphism rings, and for directing us to \cite{vig-right-adj}. 
The author was supported by the National Science Foundation under agreement DMS-2302623 and a start-up grant from Purdue University. 

\section{Some irreducible diagrams}\label{sec:irreddiag}

Let $p>3$ be a prime number. 
Let $F$ be a nonarchimedean local field with residue field $k$ a proper extension of $\F_p$ and uniformizer $\varpi$. 
We let $G = \GL_2(F)$, $K = \GL_2(\cO_F)$, and $Z \subset G$ be the center of $G$. 
Let $I\subset K$ be the upper triangular Iwahori subgroup, $I_1 \subset I$ the normal pro-$p$ Sylow subgroup, and $\Pi \in N_G(I)$ be $\left(\begin{smallmatrix}
0 & 1 \\ \varpi & 0 \end{smallmatrix} \right)$. 
For an $I$-module $V$, let $V^{(\Pi)}$ be the $\Pi$-twisted $I$-module $\{v^{(\Pi)}\mid v\in V\}$, with $I$-action given by $g\cdot v^{(\Pi)} = (\Pi g\Pi^{-1}\cdot v)^{(\Pi)}$ for $g\in I$ and $v\in V$. 

Let $E$ be a field extension of the residue field $k$ of $F$. 
Any simple $E[K]$-module is obtained via inflation from $\GL_2(k)$, is defined over $k$, and is absolutely irreducible (see e.g.~\cite[\S 1]{BL}). 
Further, the functor of taking $\GL_2(k)$-invariants, and thus also the bifunctor $\Hom_{\GL_2(k)}(-,-)$, commutes with base change from $k$ to $E$. 
In particular, formation of socle and cosocle commute with base change and indecomposability is preserved by base change. 

To construct the desired representation in the generality of Theorem \ref{thm:main}, we will use constructions (and notation) from \cite{GLS} rather than \cite{le,gs}. 
The $K$-representation $D_0$ in \cite[\S 2]{GLS} has a model over $k$ and thus a model over $E$. 
This $E[K]$-representation, which we denote $D_0$ in the present article, is rigid of Loewy length $2$ and has multiplicity free socle and cosocle. 
Let $l$ be $[k:\F_p]$ if $[k:\F_p]$ is odd and $2[k:\F_p]$ if $[k:\F_p]$ is even. 
By our assumptions, $l \geq 3$. 
The $E[K]$-socle of $D_0$ has a decomposition 
\[\sigma\oplus \sigma_1\oplus\cdots\oplus\sigma_{l-1}\oplus\sigma_1'\oplus\cdots\oplus\sigma_{l-1}'\]
into simple $E[K]$-modules. 
We also sometimes write $\sigma_0$ and $\sigma_l$ for $\sigma$. 
The $E[K]$-cosocle of $D_0$ has a decomposition 
\[\sigma^{[s]}\oplus \sigma_1^{[s]}\oplus\cdots\oplus\sigma_{l-1}^{[s]}\oplus\sigma_1^{\prime [s]}\oplus\cdots\oplus\sigma_{l-1}^{\prime [s]}\]
into simple $E[K]$-modules where for each simple $E[K]$-submodule $\kappa$ of $D_0$, $\kappa^{[s]}$ is the unique up to isomorphism simple $E[K]$-module with $(\kappa^{[s]})^{I_1}\cong (\kappa^{I_1})^{(\Pi)}$. 
Note that some Jordan-H\"older factors of $\soc_{E[K]} D_0$ may also appear in $\cosoc_{E[K]} D_0$. 

Now let $\tld{E}$ be a division algebra over $E$ with center $\tld{C}$. Eventually, we will assume that $\tld{E}$ is of countable $E$-dimension, but not at the moment. 
Let $\tld{D}_0$, $\tld{\sigma}$, etc.~denote $D_0 \otimes_E \tld{E}$, $\sigma \otimes_E \tld{E}$, etc. 
These objects have natural left and right $\tld{E}$-actions, which we will distinguish by writing on the left and right, respectively. 
The $E[K]$-socle of $\tld{D}_0$ has a decomposition 
\[\tld{\sigma}\oplus \tld{\sigma}_1\oplus\cdots\oplus\tld{\sigma}_{l-1}\oplus \tld{\sigma}_1'\oplus\cdots\oplus\tld{\sigma}_{l-1}'\]
since there is a natural isomorphism between $-\otimes_E \tld{E}$ and taking a direct sum over an $E$-basis for $\tld{E}$, and the formation of socle commutes with direct sums. 
By an abuse, we identify $D_0$, $\sigma$, etc.~with the image of the natural injection $D_0 \rightarrow \tld{D}_0$, $\sigma \rightarrow\tld{\sigma}$, etc. 
Fix a $\tld{C}$-valued character $\xi$ of $\Z \cong \left(\begin{smallmatrix}
\varpi & 0 \\ 0 & \varpi \end{smallmatrix} \right)^{\Z}$. 
We extend the induced $K$-action on $\tld{D}_0$ to $KZ \cong K \times \left(\begin{smallmatrix}
\varpi & 0 \\ 0 & \varpi \end{smallmatrix} \right)^\Z$ so that $\left(\begin{smallmatrix}
\varpi & 0 \\ 0 & \varpi \end{smallmatrix} \right)^\Z$ acts via $\xi$. 
Note that the left (and right) $\tld{E}$-action commutes with the action of $E[KZ]$. 

We let $D_1 = D_0^{I_1}$ and $\tld{D}_1 = D_1\otimes_E \tld{E} \cong \tld{D}_0^{I_1}$. 
We extend the $IZ$-action on $\tld{D}_1$ to $N_G(I)$ as follows. 
Let $S_1 = (\soc_{E[K]} D_0)^{I_1}$ and $Q_1 = (\cosoc_{E[K]} D_0)^{I_1}$ which are semisimple multiplicity-free $E[I]$-representations. 
As in the description of $\soc_{E[K]} \tld{D}_0$, we have $\tld{S}_1 = S_1\otimes_E \tld{E} \cong (\soc_{E[K]} \tld{D}_0)^{I_1}$ and $\tld{Q}_1 = Q_1\otimes_E \tld{E} \cong (\cosoc_{E[K]} \tld{D}_0)^{I_1}$. 
Moreover, there are $E[I]$-isomorphisms $D_1 \cong S_1\oplus Q_1$ and $Q_1^{(\Pi)} \cong S_1$. 
Fixing an $E[I]$-isomorphism $\psi: S_1 \cong Q_1^{(\Pi)}$ induces an $E[I]$-isomorphism $\tld{\psi}: \tld{S}_1 \cong \tld{Q}_1^{(\Pi)}$ that commutes with left and right $\tld{E}$-actions. 
Then there is a unique extension of the $IZ$-action on $\tld{D}_1$ to an $N_G(I)$-action such that $v\mapsto (\Pi v)^{(\Pi)}$ for $v\in \tld{S}_1$ is the map $\tld{\psi}$. 
Moreover, this $N_G(I)$-action commutes with both left and right $\tld{E}$-actions. 

Let $\tld{D}_{0}(\infty)=\bigoplus_{i\in\mathbb{Z}}\tld{D}_{0}(i)$ be the smooth $KZ$-representation where there is a fixed isomorphism $\tld{D}_{0}(i)\cong \tld{D}_{0}$ of $\tld{E}[KZ]$-representations for every $i\in\mathbb{Z}$. 
We denote the natural inclusion $\tld{D}_{0}\iso \tld{D}_{0}(i)\hookrightarrow \tld{D}_{0}(\infty)$ by $\iota_{i}$, and write $v_{i}=\iota_{i}(v)$ for $v\in \tld{D}_{0}$ and $i\in\mathbb{Z}$. 
Each $\iota_i$ intertwines both the left and right $\tld{E}$-actions. 
Let $\tld{D}_{1}(\infty):=\tld{D}_{0}(\infty)^{I_1}\cong\bigoplus_{i\in\mathbb{Z}}(\tld{S}_{1}\oplus \tld{Q}_{1})$. 
Fix $\lambda=(\lambda_{i})\in\prod_{i\in\mathbb{Z}} \tld{E}^{\times}$. 
We now give an extension, depending on $\lambda$, of the $IZ$-action on $\tld{D}_{1}(\infty)$ to an $N_G(I)$-action which is \emph{not} the direct sum over $\Z$ of the $N_G(I)$-action on $\tld{D}_{1}$. 
For all integers $i\in\mathbb{Z}$, let $\Pi$ act by 
\begin{equation*}
\Pi v_{i}:=
\begin{cases}
(\Pi v)_{i-1} &\text{if $v\in \tld{S}_{1}^{\chi(\sigma_{1})}$,}\\
(\Pi v)_{i+1}\lambda_i &\text{if $v\in \tld{S}_{1}^{\chi(\sigma_{1}')}$,}\\
(\Pi v)_{i} &\text{if $v\in \tld{S}_{1}^{\chi}$ for $\chi\in\{\chi(\sigma),\chi(\sigma_{2}),\ldots ,\chi(\sigma_{l-1}),\chi(\sigma_{2}'),\ldots ,\chi(\sigma_{l-1}')\}$.}
\end{cases}
\end{equation*} 
This uniquely determines a smooth $N_G(I)$-action on $\tld{D}_{1}(\infty)$ extending the action of $IZ$. 
Note that the left actions of $N_G(I)$ and $\tld{E}$ commute. 
We get a basic $0$-diagram (in the sense of \cite[\S 9]{BP} without the assumption that $\left(\begin{smallmatrix}
\varpi & 0 \\ 0 & \varpi \end{smallmatrix} \right)^\Z$ acts trivially) $\tld{D}(\lambda):=(\tld{D}_{0}(\infty),\tld{D}_{1}(\infty),\mathrm{can})$ where can is the canonical inclusion $\tld{D}_{1}(\infty)|_{IZ}\hookrightarrow \tld{D}_{0}(\infty)|_{IZ}$. 

In what follows, we use the convention that $\prod_{j=0}^{m-1} \lambda_j = \lambda_0\lambda_1\cdots \lambda_{m-1}$. 

\begin{proposition}\label{prop:irred-diagram} If $\lambda_0 \in E$, $\lambda_{i}\neq\lambda_{0}$ for all $i\neq 0$, and the $E$-span of $\{\prod_{j=0}^{m-1} \lambda_j\mid m \geq 1 \}$ is $\tld{E}$, then the basic $0$-diagram $\tld{D}(\lambda)$ is irreducible over $E$. 
\end{proposition}
\begin{proof}
The proof is a fiber product of the proofs of \cite[Proposition 3.1]{GLS} and \cite[Theorem 3.4]{le}. 
Let $W\subset \tld{D}_{0}(\infty)$ be a nonzero $E[KZ]$-subrepresentation such that $\Pi$ stabilizes $W^{I_1} \subset \tld{D}_1(\infty)$. 
We claim that $W=\tld{D}_{0}(\infty)$. 
We begin with a series of lemmata. 


\begin{lemma}\label{lemma:shift}
    If $(c_{i})\in\bigoplus_{i\in\mathbb{Z}}\tld{E}$ such that $\left(\sum_{i}\iota_{i}c_{i}\right)(\sigma)\subseteq W$, then $\left(\sum_{i}\iota_{i-m}c_{i}\right)(\sigma)\subseteq W$ and 
    \[\left(\sum_{i}\iota_{i+m}\left(c_{i}\prod_{j=0}^{m-1} \lambda_{i+j}\right)\right)(\sigma)\subseteq W\] 
    for all $m>0$. 
\end{lemma}
\begin{proof}
    Let $M(\sigma^{[s]})$ be the unique $E[K]$-direct summand of $D_0$ with cosocle isomorphic to $\sigma^{[s]}$ (see \cite[\S 2]{GLS}). 
    The $\Pi$-action on $\left(\sum_{i}\iota_{i}c_{i}\right)(S_{1}^{\chi(\sigma)})$ gives that $\left(\sum_{i}\iota_{i}c_{i}\right)(Q_{1}^{\chi(\sigma)^{s}})\subseteq W^{I_1}$ which implies that 
    \[
    \left(\sum_{i}\iota_{i}c_{i}\right)(M(\sigma^{[s]}))\subseteq W\] 
    because $Q_1^{\chi(\sigma)^s}$ generates $M(\sigma^{[s]})$ as an $E[K]$-module. 
    Hence
\begin{equation}\label{eq1}
\left(\sum_{i}\iota_{i}c_{i}\right)(\soc_{E[K]} M(\sigma^{[s]})) = \left(\sum_{i}\iota_{i}c_{i}\right)(\sigma_{1}\oplus\sigma_{1}')\subseteq W.
\end{equation} 
Now the $\Pi$-action on $\left(\sum_{i}\iota_{i}c_{i}\right)(S_{1}^{\chi(\sigma_{1})})$ and $\left(\sum_{i}\iota_{i}c_{i}\right)(S_{1}^{\chi(\sigma_{1}')})$ gives respectively \[\left(\sum_{i}\iota_{i-1}c_{i}\right)(Q_{1}^{\chi(\sigma_{1})^{s}})\subseteq W^{I_1}\hspace{2mm}\mathrm{and}\hspace{2mm}
\left(\sum_{i}\iota_{i+1}c_{i}\lambda_{i}\right)(Q_{1}^{\chi(\sigma_{1}')^{s}})\subseteq W^{I_1}.\]
(We use here that $v \cdot \lambda_i = \lambda_i \cdot v$ for $v\in Q_1^{\chi(\sigma_{1}')^{s}}$.)
Recalling from \cite[\S 2]{GLS} that $E(\sigma_{2},\sigma_{1}^{[s]}) \subset D_0$ is the unique $E[K]$-direct summand which is a nonsplit extension of $\sigma_{1}^{[s]}$ by $\sigma_{2}$, we have  
\begin{equation*}
\left(\sum_{i}\iota_{i-1}c_{i}\right)(E(\sigma_{2},\sigma_{1}^{[s]}))\subseteq W\hspace{2mm}\mathrm{and}\hspace{2mm}\left(\sum_{i}\iota_{i+1}c_{i}\lambda_{i}\right)(E(\sigma_{2}',\sigma_{1}'^{[s]}))\subseteq W 
\end{equation*} 
so that 
\begin{equation*}
\left(\sum_{i}\iota_{i-1}c_{i}\right)(\sigma_2)\subseteq W\hspace{2mm}\mathrm{and}\hspace{2mm}\left(\sum_{i}\iota_{i+1}c_{i}\lambda_{i}\right)(\sigma_2')\subseteq W. 
\end{equation*} 
Acting by $\Pi$ again and iterating a similar process shows that 
\begin{equation*}
\left(\sum_{i}\iota_{i-1}c_{i}\right)(\sigma)\subseteq W\hspace{2mm}\mathrm{and}\hspace{2mm}\left(\sum_{i}\iota_{i+1}c_{i}\lambda_{i}\right)(\sigma)\subseteq W. 
\end{equation*} 
The lemma now follows inductively. 
\end{proof}

\begin{lemma}\label{lemma:iota0}
    We have $\Hom_{E[K]}(\sigma,W\cap\iota_0(\tld{D}_0)) \neq 0$.
\end{lemma}
\begin{proof}
    We have $\Hom_{E[K]}(\tau, W)\neq 0$ for some simple $\tau\subset\mathrm{soc}_{E[K]}D_{0}$. 
If $\tau = \sigma_m$ (resp.~$\tau = \sigma_m'$) with $1\leq m<l-1$, then we have $\Hom_{E[K]}(\kappa,W) \neq 0$ for $\kappa = \sigma_{m+1}$ (resp.~$\kappa = \sigma_{m+1}'$) as in the proof of Lemma \ref{lemma:shift}. 
Similarly, if $\tau = \sigma_{l-1}$ or $\sigma_{l-1}'$, then $\Hom_{E[K]}(\sigma,W)\neq 0$. 
Without loss of generality, we can assume that $\tau = \sigma$. 
Since \[0\neq \Hom_{E[K]}(\sigma,W)\subset \Hom_{E[K]}(\sigma,\tld{D}_0(\infty)) \cong \bigoplus_{i\in\mathbb{Z}}\Hom_{E[K]}(\sigma,\iota_i (\tld{D}_0)) \cong \bigoplus_{i\in\mathbb{Z}} \iota_i|_\sigma\cdot\tld{E},\]
there exists a nonzero $(c_{i})\in\bigoplus_{i\in\mathbb{Z}}\tld{E}$ such that 
\begin{equation}\label{eqn:sigma}
    \left(\sum_{i}\iota_{i}c_{i}\right)(\sigma)\subset W.
\end{equation} 
Suppose that $\#(c_{i}):=\#\{i \in \Z\mid c_i\neq 0\}$ is minimal amongst nonzero $(c_{i})\in\bigoplus_{i\in\mathbb{Z}}\tld{E}$ satisfying \eqref{eqn:sigma}. 
By Lemma \ref{lemma:shift} we can assume without loss of generality that $c_0$ is nonzero. 
Lemma \ref{lemma:shift} implies that $\left(\sum_{i} \iota_{i+1}c_{i}\lambda_i\right)(\sigma)\subset W$ and then that $\left(\sum_{i} \iota_{i}c_{i}\lambda_i\right)(\sigma)\subset W$. 
Using that $\lambda_0 \in E^\times\subset \tld{C}^\times$, we have that $\left(\sum_{i} \iota_{i}c_{i}(\lambda_i-\lambda_0)\right)(\sigma)\subset W$. 
Since $\lambda_i-\lambda_0 \neq 0$ for all $i \neq 0$, $\#(c_{i}(\lambda_i-\lambda_0)) = \#(c_i)-1$. 
Minimality of $\#(c_i)$ implies that $c_i = 0$ for $i\neq 0$. 
\end{proof}

\begin{lemma}\label{lemma:sigma0}
    We have $\iota_0(\tld{\sigma})\subset W$. 
\end{lemma}
\begin{proof}
    We have that $\iota_0 c_0(\sigma) \subset W$ for some $c_0 \in \tld{E}^\times$ by Lemma \ref{lemma:iota0}. 
    By Lemma \ref{lemma:shift}, we also have that $\iota_0(c_0\prod_{j=0}^{m-1} \lambda_j)(\sigma) \subset W$ for all $m \geq 1$. 
    Since the $E$-span of $\{c_0\prod_{j=0}^{m-1} \lambda_j \mid m\geq 1\}$ is $\tld{E}$, we are done. 
\end{proof}
Lemmas \ref{lemma:sigma0} and \ref{lemma:shift} imply that $\iota_i(\tld{\sigma}) \subset W$ for all $i \in \Z$ so that, arguing as in the proof of Lemma \ref{lemma:shift}, $\iota_i(\tld{D}_0) \subset W$ for all $i \in \Z$. 
We conclude that $W = \tld{D}_0(\infty)$. 
\end{proof}

\section{Proof of Theorem \ref{thm:main}}

Let $\lambda = (\lambda_i)_{i\in \Z} \in \prod_{i\in \Z}\tld{E}^\times$ and $\tld{D}(\lambda) = (D_0(\infty),D_1(\infty),\mathrm{can})$ as in \S \ref{sec:irreddiag}.

\begin{theorem}\label{thm:irred} 
        \begin{enumerate}
            \item \label{item:extend} There is an $E[G]$-representation $\pi$ with a commuting left $\tld{E}$-action and an $\tld{E}$-linear injection of diagrams $\tld{D}(\lambda) = (\tld{D}_0(\infty),\tld{D}_1(\infty),\mathrm{can}) \hookrightarrow (\pi|_{KZ},\pi|_{N_G(I)},\mathrm{can})$ such that the map $\tld{D}_0(\infty) \hookrightarrow \pi|_{KZ}$ induces an isomorphism on $E[K]$-socles and $\tld{D}_0(\infty)$ generates $\pi$ as an $E[G]$-representation. 
            \item \label{item:hecke} For $\pi$ as in \eqref{item:extend} and any irreducible smooth $E[K]$-representation $\sigma'$ and $\varphi \in \Hom_{E[K]}(\sigma',\pi)$, we have $\varphi \cdot \mathcal{H}(\sigma') = E\cdot \xi\left(\begin{smallmatrix}
\varpi & 0 \\ 0 & \varpi \end{smallmatrix} \right)^{\Z}\cdot \varphi$ where $\mathcal{H}(\sigma')$ denotes the Hecke algebra $\End_{E[G]}(\mathrm{ind}_{K}^{G} \sigma')$. 
            In particular, when $\pi$ is an irreducible $E[G]$-representation, the following are equivalent: 
            \begin{itemize}
                \item $\pi$ admits a Hecke eigenvalue;
                \item $\pi$ has a central character; and
                \item $\left(\begin{smallmatrix}
\varpi & 0 \\ 0 & \varpi \end{smallmatrix} \right)$ acts on $\pi$ by a scalar in $E^\times$. 
        \end{itemize}
            \item \label{item:irred} If $\lambda$ is as in Proposition \ref{prop:irred-diagram}, then any $\pi$ as in \eqref{item:extend} is an irreducible smooth $E[G]$-representation. 
            \item \label{item:endopi} Suppose that $\lambda_0 \in E$, $\lambda_i \neq \lambda_0$ for all $i\neq 0$, and $\pi$ is as in \eqref{item:extend} and is irreducible as an $E[G]$-representation. 
            Then $\End_{E[\GL_2(F)]}(\pi) \cong \tld{E}$ as $E$-algebras. 
        \end{enumerate}  
\end{theorem}
\begin{proof}
    To prove \eqref{item:extend}, we adapt the proof of \cite[Theorem 3.2]{le}. 
    Let $D_0|_K \hookrightarrow \Omega$ denote a $E[K]$-injective envelope. 
    We first claim that $\Hom_{E[I]}(\chi,\Omega|_I) = 0$ for any character $\chi:I \rightarrow E^\times$ such that $\chi=\chi^{(\Pi)}$ (recall from \S \ref{sec:irreddiag} that $\chi^{(\Pi)}(g) = \chi(\Pi g \Pi^{-1})$). 
    By Frobenius reciprocity, we have $\Hom_{E[I]}(\chi,\Omega|_I) \cong \Hom_{E[K]}(\Ind_I^K \chi,\Omega)$. 
    Under this hypothesis of $\chi$, $\Ind_I^K \chi$ is semisimple (isomorphic to the direct sum of twists of the trivial and Steinberg simple $E[K]$-modules) so that any element of $\Hom_{E[K]}(\Ind_I^K \chi,\Omega)$ factors through $\soc_{E[K]} \Omega \cong \soc_{E[K]} D_0$. 
    Another application of Frobenius reciprocity (and the containment $I_1\subset \ker\chi$) gives $\Hom_{E[I]}(\chi,\Omega|_I) \cong \Hom_{E[I]}(\chi,D_1) = 0$ where the last equality follows from the description of $D_1$ in \cite[\S 2]{GLS}. 
    
    It is well-known (see e.g.~the proof of \cite[Lemma 9.6]{BP}) using Frobenius reciprocity that 
    \[
    \dim_E \Hom_{E[I]}(\chi,\Omega|_I) = \dim_E \Hom_{E[I]}(\chi^{(\Pi)},\Omega|_I)
    \]
    for every character $\chi: I \rightarrow E^\times$. 
    With the claim from the previous paragraph, there exist semisimple $E[I]$-submodules $S_1',Q_1' \subset \Omega|_I$ and isomorphisms $\soc_{E[I]} \Omega \cong S_1\oplus Q_1 \oplus S_1' \oplus Q_1'$ and $S_1' \cong Q_1^{\prime (\Pi)}$. \
    This gives an extension of the earlier $\psi: S_1 \cong Q_1^{(\Pi)}$ to $\psi: S_1\oplus S_1' \cong Q_1^{(\Pi)} \oplus Q_1^{\prime(\Pi)}$. 
    There is a direct sum decomposition $\Omega|_I \cong \Omega_S \oplus \Omega_Q$ where $S_1\oplus S_1' \subset \Omega_S$ and $Q_1\oplus Q_1' \subset \Omega_Q$ are $E[I]$-injective envelopes. 
    By the $E[I]$-injectivity of $\Omega_Q$, $\psi$ extends to a map $\psi: \Omega_S \rightarrow \Omega_Q^{(\Pi)}$ which is necessarily an isomorphism as it is a map between $E[I]$-injectives which is an isomorphism on $E[I]$-socles. 
    For later use, we fix $S_1^{\chi(\sigma_1)} \subset \Omega_{\chi(\sigma_1)} \subset \Omega_S \subset \Omega$ and $S_1^{\chi(\sigma_1')} \subset \Omega_{\chi(\sigma_1')} \subset \Omega_S \subset \Omega$, where the first respective inclusions are $E[I]$-injective envelopes, and we fix a direct sum decomposition $\Omega \cong \Omega_{\chi(\sigma_1)} \oplus \Omega_{\chi(\sigma_1')} \oplus \Omega'$ with $\Omega' \subset \Omega$ an $E[I]$-injective direct summand. 

    Let $\tld{\Omega}$ be $\Omega \otimes_E \tld{E}$. 
    We extend the induced $E[K]$-action on $\tld{\Omega}$ to $E[KZ]$ so that $\left(\begin{smallmatrix}
    \varpi & 0 \\ 0 & \varpi \end{smallmatrix} \right)^\Z$ acts through the character $\xi$. 
    Further, there is a unique extension of the $IZ$-action on $\tld{\Omega}$ to an $N_G(I)$-action such that $v\mapsto (\Pi v)^{(\Pi)}$ for $v\in \tld{\Omega}_S:=\Omega_S\otimes_E \tld{E}$ is the isomorphism $\tld{\psi}: \tld{\Omega}_S \cong \tld{\Omega}_Q^{(\Pi)}$ induced by $\psi$. 
    
    We now repeat for the diagram $\tld{\Omega} \subset \tld{\Omega}$ many of the constructions from \S \ref{sec:irreddiag} for the diagram $\tld{D}_1 \subset \tld{D}_0$.     
    We let $\tld{\Omega}(\infty)$ be the $\tld{E}[KZ]$-module $\oplus_{i\in \Z} \tld{\Omega}(i)$ where there is a fixed $\tld{E}[KZ]$-isomorphism $\tld{\Omega}(i) \cong \tld{\Omega}$, and as before we denote by $\iota_i$ the inclusion $\tld{\Omega} \cong \tld{\Omega}(i) \subset \tld{\Omega}(\infty)$ and write $v_i$ for $\iota_i(v)$ for $v\in \tld{\Omega}$. 
    We extend the $IZ$-action on $\tld{\Omega}(\infty)$ to an $N_G(I)$-action by the formula 
    \begin{equation*}
    \Pi v_{i}:=
    \begin{cases}
    (\Pi v)_{i-1} &\text{if $v\in \tld{\Omega}_{\chi(\sigma_1)}$,}\\
    (\Pi v)_{i+1}\lambda_i &\text{if $v\in \tld{\Omega}_{\chi(\sigma_1')}$,}\\
    (\Pi v)_{i} &\text{if $v\in \tld{\Omega}'$,}
    \end{cases}
    \end{equation*} 
    where $\tld{\Omega}_{\chi(\sigma_1)} = \Omega_{\chi(\sigma_1)}\otimes_E \tld{E}$, etc.~as usual. 
    As before, this is \emph{not} the direct sum over $\Z$ of the $N_G(I)$-action on $\tld{\Omega}$. 
    Note that the left $\tld{E}$-action commutes with the actions of $N_G(I)$ and $KZ$ on $\tld{\Omega}(\infty)$. 
    There is an obvious (left) $\tld{E}$-linear inclusion $(\tld{D}_0(\infty),\tld{D}_1(\infty),\mathrm{can}) \subset (\tld{\Omega}(\infty),\tld{\Omega}(\infty),\mathrm{id}_{\tld{\Omega}(\infty)|_{IZ}})$ of diagrams. \cite[Corollary 5.5.5]{Paskunas} gives an $E[G]$-action on $\tld{\Omega}(\infty)$ simultaneously extending the actions of $E[KZ]$ and $E[N_G(I)]$. 
    Moreover, the natural left $\tld{E}$-action commutes with the $E[G]$-action since it commutes with the actions of the subgroups $KZ$ and $N_G(I)$ which jointly generate $G$. 
    We let $\pi\subset \tld{\Omega}(\infty)$ denote the $E[G]$-subrepresentation generated by $\tld{D}_0(\infty)$. 
    We have the inclusion $(\tld{D}_0(\infty),\tld{D}_1(\infty),\mathrm{can}) \hookrightarrow (\pi|_{KZ},\pi|_{N_G(I)},\mathrm{can})$. 
    Since $\tld{D}_0(\infty)$ is stable under the left $\tld{E}$-action, so is $\pi$. 
    Thus $\pi$ has commuting left actions of $E[G]$ and $\tld{E}$ and the inclusion $\tld{D}_0(\infty) \hookrightarrow \pi$ is $\tld{E}$-linear and $E[K]$-linear. 
    It is clear that the inclusions $\tld{D}_0(\infty) \subset \pi \subset \tld{\Omega}(\infty)$ induce isomorphisms on $E[K]$-socles. 
    This completes the proof of \eqref{item:extend}. 

    Next, we prove \eqref{item:hecke}. We reduce to the case where $\varphi$ is nonzero so that $\sigma'$ is isomorphic to $\sigma_j$ or $\sigma_j'$ with $0 \leq j \leq l-1$. 
    The Hecke algebra $\mathcal{H}(\sigma')$ is generated over $E$ by the unique up to scalar Hecke operator $T$ with support equal to the $K$-double coset containing $\left(\begin{smallmatrix}
1 & 0 \\ 0 & \varpi \end{smallmatrix} \right)$ and the operators acting by $\xi\left(\begin{smallmatrix}
\varpi & 0 \\ 0 & \varpi \end{smallmatrix} \right)^{\pm1}$ \cite{herzig-satake}. 
    The operator $T$ acts on $\Hom_{E[K]}(\sigma',\pi)$ by $0$ as in the first paragraph of the proof of \cite[Lemma 4.1]{GLS}. 
    This gives the first part of \eqref{item:hecke}. 
    For the equivalences, the first bullet point implies the second and the second implies the third for any irreducible smooth $G$-representation over $E$. 
    The third bullet point implies that $\varphi \cdot \mathcal{H}(\sigma') = E \cdot \varphi$ for any $\sigma'$ and nonzero $\varphi$ by the first part of \eqref{item:hecke} which implies the first bullet point. 
    
    We now prove \eqref{item:irred}. 
    If $\pi' \subset \pi$ is a nonzero $E[G]$-subrepresentation, then $(\pi'|_{KZ} \cap \tld{D}_0(\infty),\pi'|_{KZ} \cap \tld{D}_1(\infty),\mathrm{can}) \subset (\tld{D}_0(\infty),\tld{D}_1(\infty),\mathrm{can})$ is a nonzero basic $0$-subdiagram over $E$. 
    Indeed, $\pi'|_{KZ} \cap \tld{D}_0(\infty)$ contains the $E[K]$-socle of $\pi'|_{K}$ and is thus nonzero, and it has $I_1$-invariants $\pi'|_{KZ} \cap \tld{D}_1(\infty)$. 
    Then we have that $\pi'|_{KZ}$ contains $\tld{D}_0(\infty)$ by Proposition \ref{prop:irred-diagram} so that $\pi'$ contains the $E[G]$-representation generated by $\tld{D}_0(\infty)$, which is $\pi$. 
    
    Finally, we prove \eqref{item:endopi}. 
    Assume the hypotheses therein. 
    The left action of $\tld{E}$ on $\pi$ gives an injection $\tld{E} \hookrightarrow \End_{E[\GL_2(F)]}(\pi)$ of $E$-algebras. 
    We will show that this map is an isomorphism. 
    Since $\pi$ is irreducible, the restriction map $\End_{E[\GL_2(F)]}(\pi) \rightarrow \Hom_{E[\GL_2(\cO_F)]}(\sigma,\pi)$ induced by $\iota_0: \sigma \hookrightarrow \pi$ is injective. 
    We will show that this map factors through $\Hom_{E[\GL_2(\cO_F)]}(\sigma,\iota_0(\tld{\sigma})) \cong \tld{E}$. 
    Granting this, the composition $\tld{E} \hookrightarrow \End_{E[\GL_2(F)]}(\pi) \hookrightarrow \tld{E}$ would be the identity morphism by the $\tld{E}$-linearity in \eqref{item:extend}, and we would be done.  

    We now show the above factorization. 
    Let $v^j$ be a nonzero element in $\sigma_j^{I_1}$ for each $j = 0,\ldots,l-1$ and $w^j$ be a nonzero element in $(\sigma'_j)^{I_1}$ for each $j = 1,\ldots,l-1$. By the proof of Lemma \ref{lemma:shift}, there is an element $R_j \in E[\GL_2(\cO_F)]\Pi$ (resp.~$R_j' \in E[\GL_2(\cO_F)]\Pi$) for each $j = 0,\ldots,l-1$ such that 
    \begin{itemize}
        \item $R_j \iota_i(v^j) = \iota_i(v^{j+1})$ (resp.~$R'_j \iota_i(w^j) = \iota_i(w^{j+1})$) for all $i\in \Z$ and $j=2,\ldots,l-1$ (with $v^l$ and $w^l$ taken to be $v^0$); 
        \item $R_1 \iota_i(v^1) = \iota_{i-1}(v^2)$ (resp.~$R'_1 \iota_i(w^1) = \iota_{i+1}(w^2 \lambda_i)= \iota_{i+1}(\lambda_i w^2)$) for all $i\in \Z$; and
        \item $R_0 \iota_i(v^0) = \iota_i(v^1)$ (resp.~$R'_0 \iota_i(v^0) = \iota_i(w^1)$) for all $i\in \Z$. 
    \end{itemize} 
    Let $S = R_{l-1}R_{l-2}\cdots R_0R'_{l-1}R'_{l-2}\cdots R'_0 \in E[\GL_2(F)]$. 
    Then $S \iota_i(v^0) = \iota_i(\lambda_i v^0)$ for all $i\in \Z$. 

    Let $\varphi \in \End_{E[\GL_2(F)]}(\pi)$. 
    Now $\varphi \circ \iota_0(v^0)$ is necessarily contained in the $I_1$-invariant part of the $\sigma$-isotypic part of $\soc_{E[\GL_2(\cO_F)]} \pi$, which is $\oplus_i \iota_i(\tld{\sigma}^{I_1})$, and is thus of the form $\sum_i \iota_i c_i(v^0)$ for $(c_i) \in \oplus_i \tld{E}$. 
    Then 
    \[
    S(\varphi\circ \iota_0 (v^0)) = S(\sum_i \iota_i c_i(v^0)) = \sum_i \iota_i c_i \lambda_i(v^0). 
    \]
    On the other hand, $S(\varphi\circ \iota_0 (v^0)) = \varphi(S\iota_0 (v^0)) =  \varphi ( \iota_0\lambda_0(v^0)) = \sum_i \iota_i c_i \lambda_0(v^0)$ (recall that $\lambda_0 \in E^\times\subset \tld{C}^\times$). 
    Since $\lambda_i \neq \lambda_0$ for all $i\neq 0$, we conclude that $c_i = 0$ for all $i\neq 0$. 
    Using that $v^0$ generates $\sigma$ as an $E[\GL_2(\cO_F)]$-module, this shows that the restriction map $\End_{E[\GL_2(F)]}(\pi) \rightarrow \Hom_{E[\GL_2(\cO_F)]}(\sigma,\pi)$ indeed factors through $\Hom_{E[\GL_2(\cO_F)]}(\sigma,\iota_0(\tld{\sigma}))$. 
\end{proof}

\begin{remark}
        The proof of Theorem \ref{thm:irred}\eqref{item:endopi} can be modified to show that the endomorphism ring of the diagram $\tld{D}(\lambda)$ is $\tld{E}$ when $\lambda \in \prod_i \tld{E}^\times$ with $\lambda_0 \in E^\times$ and $\lambda_i \neq \lambda_0$ for all $i\neq 0$ and $\tld{D}(\lambda)$ is an irreducible basic $0$-diagram. 
\end{remark}

\begin{proposition}\label{prop:parabolic}
    Let $F$ be a nonarchimedean local field with residue characteristic $p$. 
    Let $E$ be a field of characteristic $p$. 
Let $\ovl{P} \subset \GL_n(F)$ be a standard lower triangular parabolic subgroup with radical $\ovl{N}$ and Levi quotient $M = \ovl{P}/\ovl{N}$. 
Let $\ovl{N}^0$ be $\ovl{N} \cap \GL_n(\cO_F)$ and $M^0\subset M$ be the image of $P \cap \GL_n(\cO_F)$. 

Let $\pi$ be an $E[M]$-representation such that 
\begin{itemize}
    \item the irreducible $E[\GL_n(\cO_F)]$-subrepresentations of $\Ind_{\ovl{P}}^{\GL_n(F)} \pi$ are $M$-regular; and 
    \item denoting by 
    \[
    \mathrm{ev}_\tau: \Hom_{E[\GL_n(\cO_F)]}(\tau,\Ind_{\ovl{P}}^{\GL_n(F)}\pi) \tilde{\rightarrow} \Hom_{E[M^0]}(\tau_{\ovl{N}^0},\pi) 
    \]
    the isomorphism of right $\mathcal{H}(\tau) := \End_{E[\GL_n(F)]}(\mathrm{ind}_{\GL_n(\cO_F)}^{\GL_n(F)} \tau)$-modules induced by Frobenius reciprocity where $\mathcal{H}(\tau)$ acts on the codomain via the partial Satake homomorphism $\mathcal{S}: \mathcal{H}(\tau) \rightarrow \mathcal{H}(\tau_{\ovl{N}^0})$ (see \cite[Lemma 2.14]{herzig}),  
    we have that for any irreducible $E[\GL_n(\cO_F)]$-representation $\tau$ and $\psi \in \Hom_{E[\GL_n(\cO_F)]}(\tau,\Ind_{\ovl{P}}^{\GL_n(F)}\pi)$, $\mathrm{ev}_\tau (\psi\cdot \mathcal{H}(\tau)) = \mathrm{ev}_\tau(\psi) \cdot \mathcal{H}(\tau_{\ovl{N}^0})$. 
\end{itemize}
If $\pi$ is irreducible as an $E[M]$-representation, then $\Ind_{\ovl{P}}^{\GL_n(F)} \pi$ is irreducible as an $E[\GL_n(F)]$-representation. 
Moreover, if $\pi$ is irreducible, then $\pi$ has a Hecke eigenvalue as an $E[M]$-representation if and only if $\Ind_{\ovl{P}}^{\GL_n(F)} \pi$ has a Hecke eigenvalue as an $E[\GL_n(F)]$-representation. 
\end{proposition}
\begin{proof}
    Assuming the first claim about irreducibility of $\Ind_{\ovl{P}}^{\GL_n(F)} \pi$, the second claim about Hecke eigenvalues follows from the second bullet point. 

    To show that $\Ind_{\ovl{P}}^{\GL_n(F)} \pi$ is irreducible as an $E[\GL_n(F)]$-representation, it suffices to show that the image of any nonzero homomorphism $\psi: \tau \rightarrow \Ind_{\ovl{P}}^{\GL_n(F)} \pi$, where $\tau$ is an irreducible $E[\GL_n(\cO_F)]$-representation, generates the target as an $E[\GL_n(F)]$-representation. 
    The inverse to $\mathrm{ev}_\tau$ is given (via Frobenius reciprocity and applying $\Ind_{\ovl{P}}^{\GL_n(F)}-$) by the composition
    \begin{equation}\label{eqn:HV}
        \tau \hookrightarrow \mathrm{ind}_{\GL_n(\cO_F)}^{\GL_n(F)} \tau \hookrightarrow (\mathrm{ind}_{\GL_n(\cO_F)}^{\GL_n(F)} \tau) \otimes_{\mathcal{H}(\tau),\mathcal{S}} \mathcal{H}(\tau_{\ovl{N}^0}) \cong \Ind_{\ovl{P}}^{\GL_n(F)} \mathrm{ind}_{M^0}^{M} \tau_{\ovl{N}^0} \rightarrow \Ind_{\ovl{P}}^{\GL_n(F)} \pi
    \end{equation}
    (see \cite{HV}, particularly Theorem 4.5). 
    Suppose that \eqref{eqn:HV} is $\psi$. 
    We first claim that the image of $\mathrm{ind}_{\GL_n(\cO_F)}^{\GL_n(F)} \tau$ in $\Ind_{\ovl{P}}^{\GL_n(F)} \pi$ agrees with that of $(\mathrm{ind}_{\GL_n(\cO_F)}^{\GL_n(F)} \tau) \otimes_{\mathcal{H}(\tau)}\mathcal{H}(\tau_{\ovl{N}^0})$. 
    Indeed, $(\mathrm{ind}_{\GL_n(\cO_F)}^{\GL_n(F)} \tau) \otimes_{\mathcal{H}(\tau),\mathcal{S}} \mathcal{H}(\tau_{\ovl{N}^0})$ is generated as an $E[\GL_n(F)]$-representation by the subspace $\tau \otimes_{\mathcal{H}(\tau),\mathcal{S}} \mathcal{H}(\tau_{\ovl{N}^0})$ whose image in $\Ind_{\ovl{P}}^{\GL_n(F)} \pi$ is the image of the evaluation map $(\psi\cdot \mathcal{H}(\tau_{\ovl{N}^0})) \otimes_E \tau \rightarrow \Ind_{\ovl{P}}^{\GL_n(F)} \pi$ where $\mathcal{H}(\tau_{\ovl{N}^0})$ acts on $\psi$ via $\mathrm{ev}_\tau$. 
    By the second bullet point in the statement of the proposition, this is the image of $(\psi\cdot \mathcal{H}(\tau)) \otimes_E \tau \rightarrow \Ind_{\ovl{P}}^{\GL_n(F)} \pi$ which is contained in the image of $\mathrm{ind}_{\GL_n(\cO_F)}^{\GL_n(F)} \tau$. 
    This establishes the claim. 
    Since $\Ind_{\ovl{P}}^{\GL_n(F)} -$ is exact and $\pi$ is irreducible, the last map of \eqref{eqn:HV} is surjective, and we are done. 
    \end{proof}

\begin{theorem}\label{thm:endo}
    Let $n>1$ and $p>3$ be a prime. 
    Let $F$ be a nonarchimedean local field with residue field a proper extension of $\F_p$. 
    Let $E$ be an extension of the residue field of $F$ and $\tld{E}$ be a division algebra over  $E$ of countable dimension with center $\tld{C}$. 
    Then there exists a nonadmissible irreducible smooth $\GL_n(F)$-representation over $E$ whose endomorphism algebra is isomorphic to $\tld{E}$ as $E$-algebras. 
    Moreover, 
    \begin{itemize}
        \item there exists such a representation with a Hecke eigenvalue (and thus a central character); 
        \item if $E \subsetneq \tld{C}$, there exists such a representation without a central character; and 
        \item if $n>3$ and $E \subsetneq \tld{C}$, a representation with a central character and without a Hecke eigenvalue. 
    \end{itemize} 
\end{theorem}
\begin{proof}
    We can choose 
    \[(\lambda_i)_{i\in \Z} \in \prod_{i\in \Z} \tld{E}^\times\] 
    so that $\lambda_0 \in E$, $\lambda_{i}\neq\lambda_{0}$ for all $i\neq 0$, and the $E$-span of $\{\prod_{j=0}^{m-1} \lambda_j\mid m \geq 1 \}$ is $\tld{E}$. 
    For example, we can choose an $E$-spanning sequence $(\mu_i)_{i\geq 0}$ for $\tld{E}$ with $\mu_0 = 1$ and $\mu_i \neq \mu_{i+1}$ for all $i\geq 0$ and set $\lambda_0 = 1$, $\lambda_i = \mu_{i-1}^{-1}\mu_i$ for all $i>0$, and $\lambda_i = -1$ for all $i<0$. 
    Let $\pi$ be as in Theorem \ref{thm:irred}. 
    Then $\pi$ is an irreducible smooth $E[G]$-representation. 
    Let $\ovl{P}$ denote the lower block triangular parabolic subgroup of $\GL_n$ with blocks of size $2,1,\ldots,1$ with Levi quotient $M$. 
    If $\chi: (F^\times)^{n-2} \rightarrow \tld{C}^\times$ is any character, then $\pi \otimes_{\tld{C}} \chi$ is an irreducible $E[M]$-representation since the restriction to the factor $\GL_2(F)$ is isomorphic to $\pi$ and thus irreducible over $E$. 
    If $\chi|_{(\cO_F^\times)^{n-2}} \cong F(a) \otimes_{\tld{C}} F(b) \otimes_{\tld{C}} F(a) \otimes_{\tld{C}}\cdots$ where $0 \leq a,b< p^f-1$ are as in \cite[Lemma 4.2]{GLS} and $F(a):\cO_F^\times \rightarrow \tld{C}^\times$ maps $\alpha \mapsto \ovl{\alpha}^a$ (with $\ovl{\alpha}\in k^\times \subset \tld{C}^\times$ denoting the reduction of $\alpha$), then the proof of \cite[Lemma 4.2]{GLS} shows that $\rho = \Ind_{\ovl{P}}^{\GL_n(F)} \pi \otimes_{\tld{C}} \chi$ contains only $M$-regular Serre weights. 
    From now on, we take $\chi|_{(\cO_F^\times)^{n-2}}$ to have this form. 
    First assume moreover that $\chi$ is $E^\times$-valued. 
    Then $\rho$ is a (nonadmissible) irreducible smooth $E[\GL_n(F)]$-representation by Proposition \ref{prop:parabolic}. 
    Indeed, for any irreducible $E[\GL_n(\cO_F)]$-representation $\tau$ and $\psi \in \Hom_{E[\GL_n(\cO_F)]}(\tau,\Ind_{\ovl{P}}^{\GL_n(F)}\pi \otimes_{\tld{C}} \chi)$, we have $\mathrm{ev}_\tau(\psi) \cdot \mathcal{H}(\tau_{\ovl{N}^0}) = E \cdot \xi\left(\begin{smallmatrix}
\varpi & 0 \\ 0 & \varpi \end{smallmatrix} \right)^{\Z} \cdot \mathrm{ev}_\tau(\psi)$ by Theorem \ref{thm:irred}\eqref{item:hecke}, and one can show that $\mathrm{ev}_\tau (\psi \cdot \mathcal{H}(\tau)) = E \cdot \xi\left(\begin{smallmatrix}
\varpi & 0 \\ 0 & \varpi \end{smallmatrix} \right)^{\Z} \cdot \mathrm{ev}_\tau(\psi)$ using \cite[Lemma 2.14]{herzig} and the explicit formula for the image of the partial Satake homomorphism $\mathcal{S}$ in the proof of \cite[Proposition 2.12]{herzig}. 
    Moreover, one can arrange that $\pi$ has a Hecke eigenvalue as an $E[\GL_2(F)]$-representation by choosing $\xi\left(\begin{smallmatrix}
\varpi & 0 \\ 0 & \varpi \end{smallmatrix} \right)$ to be in $E^\times$  again by Theorem \ref{thm:irred}\eqref{item:hecke}. 
    Then $\pi\otimes_{\tld{C}} \chi$ and thus $\rho$ has a Hecke eigenvalue by the second part of Proposition \ref{prop:parabolic}. 
    If $E\subsetneq \tld{C}$, we can instead choose $\xi\left(\begin{smallmatrix}
\varpi & 0 \\ 0 & \varpi \end{smallmatrix} \right) \in \tld{C}^\times \setminus E^\times$ so that $\varpi I_n$ does not act by a scalar in $E^\times$ and, in particular, $\rho$ does not have a central character. 

    Suppose now that $n>3$ and $E\subsetneq \tld{C}$. 
    We now choose $\xi\left(\begin{smallmatrix}
\varpi & 0 \\ 0 & \varpi \end{smallmatrix} \right)$ to be an element $X \in \tld{C}^\times \setminus E^\times$, and we choose
    \[
    \chi = \bigotimes_{j=1}^{n-2}{\chi_j}
    \]
    so that $\chi_1(\varpi) = X^{-2}$, $\chi_2(\varpi) = X$, and $\chi_j(\varpi) \in E^\times$ for all $2<j\leq n-2$. 
    Note that $\pi\otimes_{\tld{C}} \chi$ does not have a central character and thus does not admit a Hecke eigenvalue. 
    One can check, again using Theorem \ref{thm:irred}\eqref{item:hecke} and \cite[Proposition 2.12, Lemma 2.14]{herzig}, that in this case, for any irreducible $E[\GL_n(\cO_F)]$-representation $\tau$ and for any $\psi \in \Hom_{E[\GL_n(\cO_F)]}(\tau,\Ind_{\ovl{P}}^{\GL_n(F)}(\pi \otimes_{\tld{C}} \chi))$, we have $\mathrm{ev}_\tau (\psi \cdot \mathcal{H}(\tau)) = E[X,X^{-1}] \cdot \mathrm{ev}_\tau(\psi) = \mathrm{ev}_\tau(\psi) \cdot \mathcal{H}(\tau_{\ovl{N}^0})$ so that $\rho$ is a (nonadmissible) irreducible smooth $E[\GL_n(F)]$-representation by Proposition \ref{prop:parabolic}. 
    Furthermore, $\rho$ has a central character since $\varpi I_n$ acts by a scalar in $E^\times$, but does not have a Hecke eigenvalue by the last part of Proposition \ref{prop:parabolic}. 

    We finally show that $\End_{E[\GL_n(F)]}(\rho) \cong \tld{E}$ as $E$-algebras with any of the above choices of $\chi$. 
    Indeed, we have $\End_{E[\GL_n(F)]}(\rho) \cong \End_{E[M]}(\pi\otimes_{\tld{C}} \chi) \cong \End_{E[\GL_2(F)]}(\pi) \cong \tld{E}$ where the first isomorphism follows from the full faithfulness of parabolic induction \cite[Theorem 5.3(1)]{vig-right-adj} and the last isomorphism follows from Theorem \ref{thm:irred}\eqref{item:endopi}. 
\end{proof}

\begin{proof}[Proof of Theorem \ref{thm:main}]
    Theorem \ref{thm:main}\eqref{item:endo} follows immediately from Theorem \ref{thm:endo}. 
    To prove Theorem \ref{thm:main}\eqref{item:central}, we need only find a proper field extension $\tld{C} = \tld{E}$ of $E$ which is of countable degree. 
    For example, one could take $\tld{E}$ to be the field of rational functions in one variable over $E$. 
\end{proof}

\bibliographystyle{amsalpha}

\bibliography{nonadm}

\end{document}